\documentclass[a4paper,11pt]{amsart}
\usepackage{amsmath}
\usepackage{graphicx}
\usepackage{epstopdf}
\usepackage{amssymb,amscd,array}
\usepackage{color}
\usepackage{float}

\newcounter{tablegroup}
\newcounter{subtable}[tablegroup]

\newtheorem{thm}{Theorem}[section]

\newtheorem{lem}[thm]{Lemma}

\newtheorem{rem}[thm]{\bf Remark}

\numberwithin{equation}{section}

\begin{document}
\title[Topological sequence entropy on trees]
{Topological sequence entropy and topological dynamics of tree maps}

\author{Aymen Daghar and Jose S. C\'anovas}

\address{Aymen Daghar, \ University of Carthage, Faculty
	of Science of Bizerte, (UR17ES21), ``Dynamical systems and their applications'', Jarzouna, 7021, Bizerte, Tunisia.}
\email{aymendaghar@gmail.com}
\address{Jose S. C\'anovas, \ Department of Applied Mathematics and Statistics, Technical University of
Cartagena, C/ Doctor Flemming sn, 30.202, Cartagena, Spain.}
\email{jose.canovas@upct.es}

\begin{abstract}
We prove that a zero topological entropy continuous tree map always displays zero topological sequence entropy when it is restricted to its non-wandering and chain recurrent sets. In addition, we show that a similar result is not possible when the phase space is a dendrite even when we consider only the restriction on the set of periodic points.
\end{abstract}
\maketitle

\section{Introduction}

Topological sequence entropy \cite{goodman} is a conjugacy invariant that is useful for distinguishing among continuous maps with zero topological entropy. Thus, it has been used to study zero entropy systems as substitution shifts \cite{deeking,leman} and one-dimensional continuous maps \cite{canovas,smifra,hric,hric2}.

For continuous maps of the interval, positive topological sequence entropy characterizes the existence of chaos in the sense of Li and Yorke \cite{smifra}. A similar result is obtained for continuous circle maps \cite{hric2}. In addition, for both interval and circle cases, zero topological entropy maps have either $\log 2$ or zero as the supremum of their topological sequence entropies \cite{canovas,canovas3}. Similar results are obtained for continuous maps on a finite tree \cite{tan}.

For a continuous interval or circle map $f:X\rightarrow X$, it is known that their topological sequence entropy can be computed by taking the restriction of $f$ to the subset $Y=\cap_{n\geq 0}f^n(X)$ \cite{bacaji2}. However, for interval maps with zero topological entropy, we find that when we restrict our map to the chain recurrent set, and hence to the non-wandering set, its topological sequence entropy is zero \cite{canovas2}. A similar result was recently proved for continuous circle maps in \cite{kuang}, but in this case, the authors proved their result just on the non-wandering set.

This paper aims to prove the same result for continuous tree maps. More precisely, we denote by $h_{\infty}(f)$ the supremum of the topological sequence entropies of a continuous map defined on a compact metric space. In addition, $\omega(f)$ denotes the $\omega$-limit set of $f$, $\Omega (f)$ denotes its non-wandering set and $CR(f)$ denotes the set of chain recurrent points. Note that, by \cite[Theorem 3]{May}, in the case of graph maps, in particular tree maps, we have $\omega(f)=\bigcap_{n\geq 0}f^{n}(\Omega(f))$. Clearly $\omega(f)$ is a closed subset of $X$. The main aim of this paper is to prove the following result.

\begin{thm}\label{Main 2}
Let $X$ be a tree and let $f:X\rightarrow X$ be continuous with zero topological entropy. Then, $$h_{\infty}(f|_{ \omega(f)})=h_{\infty}(f|_{\Omega(f)})=h_{\infty}(f|_{CR(f)})=0.$$
\end{thm}

Although the proof of this result can be done following the outline of the those for interval maps from \cite{canovas2}, here we present an alternative one. Additionally, we construct a counterexample proving that a similar result cannot be obtained when continuous zero topological entropy maps on dendrites are considered, even when we restrict the map to the set of periodic points.

The paper is organized as follows. The next section introduces the basic definitions and notations to follow the proofs. Section \ref{sec4} is devoted to prove Theorem \ref{Main 2}. In the last section, we construct the counterexample on dendrite maps.

\section{Preliminaries}

Thorough this paper $X$ will be a compact metric space with metric $d$ and $f: X\longrightarrow X$ will be a continuous map. The pair $(X,f)$ is a \textit{dynamical system}. Let $\mathbb{Z},\ \mathbb{Z}_{+}$ and $\mathbb{N}$ be the sets of integers, non-negative integers and positive integers, respectively. For $n\in \mathbb{Z_{+}}$ denote by  $f^{n}$ the $n$-$\textrm{th}$ iterate of $f$; that is, $f^{0}=\textrm{identity}$ and $f^{n}=f\circ f^{n-1}$ if $n\in \mathbb{N}$. For any $x\in X$, the subset
$\textrm{Orb}_{f}(x) = \{f^{n}(x): n\in\mathbb{Z}_{+}\}$ is called the \textit{orbit} of $x$ (under $f$). A subset $Y\subset X$ is called \textit{$f-$invariant} (resp. strongly $f-$invariant) if $f(Y)\subset Y$ (resp., $f(Y)=Y$); it is further called a \textit{minimal set} (under $f$) if it is closed, non-empty and does not contain any $f-$invariant, closed proper non-empty subset of $X$. We define the \textit{$\omega$-limit} set of a point $x$ to be the set:
\begin{align*}
\omega_{f}(x)  = \{y\in X: \liminf_{n\to +\infty} d(f^{n}(x), y) = 0\} = \underset{n\in \mathbb{N}}\cap\overline{\{f^{k}(x): k\geq n\}}.
\end{align*}
For a dynamical system $(X,f)$, the $\omega$-limit set of $f$ is denoted by: $$\omega(f)=\bigcup_{x\in X}\omega_{f}(x).$$\\
Let us recall some notions of recurrence in dynamics (see e.g. \cite{blockcoppel}).
\begin{itemize}
\item A point $x$ in $X$ is \emph{periodic} if $f^n (x)=x$ for some $n\in \mathbb{N}$. When $n=1$ the point $x$ is said to be fixed. By $P(f)$ we denote the set of periodic points.

\item A point $x$ in $X$ is \emph{wandering} if there exists a neighborhood $U$ of $x$ (called wandering neighborhood) such that $f^{-n}(U)\cap U=\emptyset$, for every $n\in\mathbb{N}$. Otherwise, the point $x$ is said to be \emph{non-wandering}. By $\Omega(f)$ we denote the set of non-wandering points.
\item For $\epsilon>0$, an $\epsilon$-chain from $x$ to $y$ under $f$ is a finite sequence $(x_{i})_{0\leq i\leq n}$ such that $x_{0}=x$ and $x_{n}=y$ and $d(f(x_{i-1}),x_{i})<\epsilon$ for any $1\leq i\leq m$. We denote by $CR(x,f)$ the set of points of $X$ for which for any $\epsilon>0$ there exist an $\epsilon$-chain under $f$ from $x$ to $y$.
\item A point $x\in X$ is \emph{chain recurrent} if $x\in CR(x,f)$. By $CR(f)$ we denote the set of chain recurrent points.
 \end{itemize}
Recall that we have the following inclusion which are strict in general: $$\omega(f)\subset \Omega(f)\subset CR(f).$$

Following \cite[page 99]{blockcoppel}, a closed set $A$ is said to be  \emph{stable} if for each open set $U$ with $A\subset U$, there exists an open set $V$ such that $A\subset V$ and for any $x\in V,\; \mathrm{Orb}_{f}(x)\subset U$. Moreover, $A$ is said to be \emph{asymptotically stable}, if in addition there exists some open set $U_{0}$ such that $A\subset U_{0}$ and for any $x\in U_{0},\;\omega_{f}(x)\subset A$. For each $x\in X$, we denote by $Q(f,x)$ the intersection of all asymptotically stable sets which contains $\omega_{f}(x)$. Note that $x\in CR(f)$ if and only if $x\in Q(f,x)$ by \cite[Proposition 39, page 113]{blockcoppel}.

Given $\mathcal{U}$ and $\mathcal{V}$ two open covers of $X$, we define the join cover $\mathcal{U}\vee \mathcal{V}:=\{ U\cap V:U\in \mathcal{U}$ and $V\in \mathcal{V} \}$. By $\mathcal{N}(\mathcal{U})$ we denote the cardinality of a finite minimal open subcover of $\mathcal{U}$.
Let $A=(a_{i})_{i\geq 0}$ be an increasing sequence of positive integers. Following \cite{goodman}, the sequence topological entropy of $f$ relative to $A$ and $\mathcal{U}$ is defined by
$$ h_{A}(f,\mathcal{U}):=\limsup_{n\rightarrow +\infty}\frac{1}{n} \log \mathcal{N} \left( \bigvee_{i=1}^nf^{-a_{i}}(\mathcal{U}) \right).$$
The {\it topological sequence entropy} of $f$ with respect to $A$ is defined by
$$h_{A}(f):=\sup \{h_{A}(f,\mathcal{U}),\; \mathcal{U}\; \rm{open\; cover\; of\; X} \} .$$
Note that when $A=(n-1)$ we receive the standard {\it topological entropy}, denoted $h(f)$. We point out that there are equivalent definitions based on Bowen's notions of spanning and separated sets (see \cite{bowen,goodman}). Following \cite{canovas}, we denote by $h_{\infty}(f)$ the supremum of $h_{A}(f)$ over all the increasing sequence $A$ of integers.

Given two discrete dynamical systems $(X,f)$ and $(Y,g)$, we say that $g$ is a factor of $f$ if there exists a continuous onto map $\pi :X\rightarrow $ such that $\pi \circ f=g\circ \pi$. If $\pi$ is a homeomorphism we say that $f$ and $g$ are conjugated. Then $h_{A}(f)\geq h_{A}(g)$ when $\pi$ is onto and $h_{A}(f)= h_{A}(g)$ when $\pi$ is an homeomorphism.

For $Y\subset X$, let $\Delta _Y:=\{ (y,y):y\in Y\}$. Next, we introduce some notation on dynamical pairs (see \cite{Wo,null}).
\begin{itemize}
\item A pair $(x,y)\in X\times X$ is \textit{proximal} if $\liminf _{n\rightarrow \infty} d(f^n(x),f^n(y))=0$. Otherwise, the pair $(x,y)$ is said to be \textit{distal}. The pair $(x,y)$ is called \textit{asymptotic} when $\limsup _{n\rightarrow \infty} d(f^n(x),f^n(y))=0$. A pair $(x,y)$ is called a \textit{Li-Yorke pair} if it is proximal but not asymptotic. The dynamical system $(X,f)$ is called \textit{distal} if it has no proximal pairs. A subset $S$ of X with at least two points is said to be a scrambled set (of $f$) if any proper pair $(a, b) \in S \times S$ is a Li-Yorke pair. A continuous map $f : X \to X$ is called Li-Yorke chaotic if it has an uncountable scrambled set.
Denote by $P(X,f)$ the set of proximal pairs and by $A(X,f)$ the set of asymptotic pairs. Finally, $f$ is Lyapunov stable if it has equicontinuous powers.
\item Given $A_1,A_2,...,A_k\subset X$, we say that $\mathcal{I}\subset \mathbb{N}$ is an \textit{independence set} of $A_1,A_2,...,A_k$ if for any non-empty finite subset $\mathcal{J} \subset \mathcal{I}$ with $m$ elements and any $s_1,...,s_m\in \{ 1,2,...,k\}$ we have that $\cap _{j\in \mathcal{J}}f^{-j}(A_{s_j})\neq \emptyset$. Following \cite{gra}, a pair $(x,y)\in X\times X$ is called a \textit{IN-pair}  if for any neighborhoods $U,V$ of $x,y$, respectively, the pair $\{U,V\}$ has arbitrarily large finite independence set for the product map $f\times f$. Denote by $IN(X,f)$ the set of IN-pair of $(X,f)$.
\end{itemize}
It is known, see \cite{gra}, that $h_{\infty}(f)>0$ if and only if $IN(X,f)\setminus \Delta_{X}\neq \emptyset$.

Finally, we introduce some basic notions on continuum theory (see \cite{nadler}).  A \emph{continuum} is a compact connected metric space.
An \emph{arc} $I$ (resp. a \emph{circle}) is any space homeomorphic to the compact interval $[0, 1]$ (resp. to the unit circle $\mathbb{S}^{1} =\{z\in \mathbb{C}: \ \vert z\vert = 1\}$). A space is called \textit{degenerate} if it is a single point, otherwise; it is \textit{non-degenerate}. By a \textit{graph} $G$, we mean a continuum which can be written as the union of finitely many arcs such that any two arcs are either disjoint or intersect only in one or both of their endpoints.
By a \textit{tree} $T$, we mean a graph containing no simple closed curve. Given $x\in T$, we denote its valence, denoted $val(x)$, as the number of connected components of $T\setminus \{x\}$. Then, the set of branching points of $T$ is $B(T):=\{ x\in T: val(x)>2 \}$. The endpoints of $T$ is then $End(T):= \{ x\in T:val(x)=1\}$. Clearly, trees are metrizable compact spaces.

A \textit{dendrite} is a locally connected continuum which contains no simple closed curve. A \textit{dendroid} is an arc-wise connected continuum which contains no simple closed curve. A map $\phi: X\to Y$ is said to be \textit{pointwise monotone} if for any $y\in Y$ the set $\phi^{-1}(y)$ is a connected subset of $X$ and \textit{monotone} if for any connected subset $C$ of $Y\; \phi^{-1}(C)$ is a connected subset of $X$. Notice that if $\phi$ is onto and closed, then $\phi$ is monotone if and only if it is pointwise monotone (see \cite[Theorem 9, page 131]{Kura}). Note that if $X$ is a tree and $\phi: X\to Y$ is a monotone onto continuous map then $Y$ is also a tree.

\section{\label{sec4}The sequence topological entropy of tree maps}
The aim of this section is to prove Theorem \ref{Main 2}. For that  we will need the following lemmas.

\begin{lem}\label{crf} \rm{\cite[Corollary 40, Propositions 44 and 48]{blockcoppel}}
Let $(X,f)$ be a dynamical system, the following hold :
\begin{enumerate}
 \item[(i)] $f(CR(f))=CR(f)$
 \item[(ii)] For any $n\geq 0,\;  CR(f^{n})=CR(f)$
\item[(iii)] If $g=f|_{CR(f)}$, then $CR(g)=CR(f)$.
\end{enumerate}
\end{lem}

\begin{lem}\label{Deomp}
Let $X$ be a tree and let $f:X\rightarrow X$ be continuous with zero topological entropy. Then, the following hold:
\begin{enumerate}
\item[(i)] Any infinite $\omega$-limit set is contained in the orbit of a periodic arc.
\item[(ii)] For any point $a$ of an infinite $\omega$-limit set, the intersection of all periodic arc containing $a$ is either $\{a\}$ or an arc $[a,b]$, with $f^{n}([a,b])\cap [a,b]=\emptyset$, for any $n\geq 1$.
\item[(iii)] For infinitely many point $a$ of an infinite $\omega$-limit set, the intersection of all periodic arc containing $a$ is reduced to $\{a\}$.
\end{enumerate}
\end{lem}
\begin{proof}
The proof of (i) is  in \cite[Corollary 1]{Gh}. The proof of (ii) follows from (i) and \cite[Theorem 3.6]{Jian Li}. For (iii), let $a\in \omega_{f}(x)$, where $\omega_{f}(x)$ is infinite and the intersection of all periodic arc containing $a$ is not reduced to $ \{a\}$. By $(i)$, we may find an $f^{m}$ periodic arc $I$ so that $a\in I$. Thus, there exists $b\in I$ such that the intersection of all periodic arc containing $a$ is $[a,b]$. Clearly, the set of pair of points $a,b$ satisfying the above condition is countable (see also \cite[Proposition 4.11]{Jian Li}). Since, by \cite{un}, $\omega_{f}(x)$ is uncountable, we conclude that infinitely many points of $I\cap \omega_{f}(x)$ satisfies condition (iii), and the proof finishes.
\end{proof}

\begin{lem}\label{x not pf}\cite[Theorem A, (5)]{Eq tree map}
Let $X$ be a tree and let $f:X\rightarrow X$ be continuous with zero topological entropy. Then, for any $x\in CR(f)\setminus P(f)$ we have that $\omega_{f}(x)$ is infinite.
\end{lem}

\begin{lem}\label{Factor CRF}
Let $(X,f)$ be a dynamical system and let $(Y,g)$ be  a factor of $(X,f)$ thought the factor map $\pi$. Then, if  $x\in CR(f)$, then $\pi(x)\in CR(g)$.
\end{lem}
\begin{proof}
  Let $a\in CR(f)$. Since $\pi $ is uniformly continuous, we can fix $\epsilon>0$ and $\eta>0$ so that for any $x,y\in X$, if $d(x,y)<\eta$ then $d(\pi(x),\pi(y))<\epsilon$.
  As $a\in CR(f)$, we can find an $\eta$-chain $\{a=a_{0},a_{1},\dots, a_{n-1},a=a_{n}\}$ such that for any $1\leq i\leq n$ we have that $d(f(a_{i-1}),a_{i})<\eta$. Thus, we have that  $d(\pi(f(a_{i-1})),\pi(a_{i}))<\epsilon$. Then $\{\pi(a)=\pi( a_{0}),\pi( a_{1}),\dots, \pi( a_{n-1}),\pi( a)=\pi( a_{n})\}$ is an $\epsilon$-chain and so $\pi(a)\in CR(g)$.
  \end{proof}

\begin{lem}\label{xw}
Let $X$ be a tree and let $f:X\rightarrow X$ be continuous with zero topological entropy. Let $x\in CR(f)\setminus P(f)$ and let $I$ be a periodic arc such that $\omega_{f}(x)\subset \mathrm{Orb}_{f}(I)$. Then, the following statements hold:
\begin{enumerate}
\item[(i)]  $x\in \mathrm{Orb}_{f}(I)$.
\item[(ii)] For any $n\geq 0,\;  f^{-n}(x)\cap CR(f) \subset  \mathrm{Orb}_{f}(I)$.
\end{enumerate}
\end{lem}
\begin{proof}
\textit{(i)} Assume this is not the case and let $I$ be a periodic arc such that  $\omega_{f}(x)\subset \mathrm{Orb}_{f}(I)$, we have $x\notin \mathrm{Orb}_{f}(I)$. Consider $\pi$ the quotient map collapsing each element of $\mathrm{Orb}_{f}(I)$ to a point, clearly $\pi$ is a monotone map inducing a factor system on the tree $\pi(X)$ with zero topological entropy, moreover by Lemma \ref{Factor CRF} $\pi(x)\in CR(\pi(f))$. Clearly $\pi(x)\notin P(f)$ and $\omega_{\pi(f)}(\pi(x))=\pi(\mathrm{Orb}_{f}(I))$ is a periodic orbit which contradict Lemma \ref{x not pf}.

\textit{(ii)} Let $z\in [f^{-n}(x)\cap CR(f) ]$, clearly $\omega_{f}(z)=\omega_{f}(x) \subset \mathrm{Orb}_{f}(I)$, thus by (i), we have $z\in \mathrm{Orb}_{f}(I)$.
\end{proof}

\begin{lem}\label{Equ}
Let $f:X\to X$ be a tree map. Then $h(f)=0$ if and only if the map $f|_{CR(f)}$ has equicontinuous powers at any periodic point.
\end{lem}
\begin{proof}
The only if part is an immediate consequence from \cite[Theorem A, (4)]{Eq tree map}.

For the if part, let $x\in P(f)$. We may assume that $x$ is a fixed point. Suppose that $f|_{CR(f)}$ does not have equicontinuous powers at $x$. Then we can find a sequence $(x_{n})_{n\geq 0}$ of $CR(f)$ converging to $x$ and a sequence $(m_{n})_{n\geq 0}$ of positive integer such that the sequence $(f^{m_{n}}(x_{n}))_{n\geq 0}$ converges to $z \in CR(f)$ with $z\neq x$. By \cite[Theorem A, (4)]{Eq tree map}, for $n$ large enough, $x_{n}\notin P(f)$ and, by \cite[Theorem 1, (3)]{Eq tree map}, the sequence $(\omega_{f}(x_{n}))_{n\geq 0}$ converges to $\omega_{f}(x)=\{x\}$. Thus $z\notin P(f)$ because otherwise the sequence $(\omega_{f}(x_{n}))_{n\geq 0}$ would also converge to $\omega_{f}(z)$ and so $\omega_{f}(z)=\{a\}$ and $a\neq z$, which contradicts Lemma \ref{x not pf}. As $z$ is not a periodic point, we conclude that $\omega_{f}(z)$ is infinite.
On the other hand, by Lemma \ref{Deomp}, (1) we can find a periodic arc such that the orbit of this arc contains $\omega_{f}(z)$. By Lemma \ref{Deomp} (3), since $\omega_{f}(z)$ is infinite, we can pick $a\in \omega_{f}(z)\setminus B(X)$ such that $I$, the periodic arc containing $a$, does not contains any branch points and $x\notin I$. Observe that $I\setminus End(I)$ is an open subset of $X$. Clearly $I\setminus End(I)$ contains $f^{k}(z)$ for some $k\geq 0$. Thus, for $n$ large enough the sequence $f^{m_{n}+k}(x_{n})$ is contained in $I$. We conclude that for any $q\geq 0,\;   f^{qm}(f^{m_{n}+k}(x_{n})) \in I$ for $n$ large enough, where $m$ is the period of $I$. Therefore $\omega_{f^{m}}(f^{m_{n}+k}(x_{n}))\subset  I$ for $n$ large enough. Recall that $\omega_{f^{m}}(f^{m_{n}+k}(x_{n}))\subset \omega_{f}(x_{n})$ and the sequence $(\omega_{f}(x_{n}))_{n\geq 0}$ converges to $\omega_{f}(x)=\{x\}$. We conclude that $x\in I$, which is a contradiction.
\end{proof}


\begin{lem}\label{Factor SE(X,f)}\cite[Lemma 2.4, 2.5]{gra}
Let $(X,f)$ be a dynamical system. Then, the following hold:
\begin{enumerate}
\item[(i)] For any $k\geq 0,\; IN(X,f)=IN(X,f^{k})$, moreover $IN(X,f)$ is closed.
\item[(ii)] Let $(Y,g)$ be a factor of $(X,f)$ thought the factor map $\pi$. If $(x,y)\in IN(X,f)$, then $(\pi(x),\pi(y))\in IN(Y,g)$.
\end{enumerate}
\end{lem}

\begin{lem}\label{ARC} \cite{canovas2}
Let $X$ be a closed interval and let $f:X\rightarrow X$ be continuous with zero topological entropy. Then $$h_{\infty}(f|_{\omega(f)})=h_{\infty}(f|_{\Omega(f)})=h_{\infty}(f|_{CR(f)})=0.$$
\end{lem}

\subsection*{Proof of Theorem \ref{Main 2}.}

As $\omega(f)\subset \Omega(f)\subset CR(f)$, it suffices to prove that $h_{\infty}(f|_{CR(f)})=0$. The proof is given in several steps and the outline is as follows: we will assume that $h_{\infty}(f|_{CR(f)})>0$ to get a contradiction at the end.

So, we assume that $h_{\infty}(f|_{CR(f)})>0$.  Then, there is a pair $(x,y)\in IN(CR(f),f)$ such that $x\neq y$. By \cite[Lemma 5.9]{gra}, since $IN(CR(f),f)\subset IN(X,f)$, we have that $x,y\in \omega(f)$. First, we prove the next two claims.
\medskip

\textit{Claim 1} At least $\omega_{f}(x)$ or $\omega_{f}(y)$ is infinite.

Suppose this is not the case. Then, by Lemma \ref{x not pf}, $x$ and $y$ are periodic points. Note that we may assume that they are fixed points. Since $(x,y)\in IN(CR(f),f)$, we can find a sequence $(x_{n})_{n\geq 0}$ of $CR(f)$ and two sequences of positive integers $(i_{n})_{n\geq 0},\;  (j_{n})_{n\geq 0}$ such that for any $n\geq 0,\; i_{n}<j_{n}$ and $f^{i_{n}}(x_{n})$ (resp. $f^{j_{n}}(x_{n})$) converges to $x$ (resp. $y$). By Lemma \ref{Equ}, $f|_{CR(f)}$ has equicontinuous powers at $x$. Then, the sequence $(f^{j_{n}-i{n}}(f^{i_{n}}(x_{n})))_{n\geq 0}$ converges to $x$, thus $x=y$ which is a contradiction. This finish the proof of the claim and so we can assume that  $\omega_{f}(x)$ is infinite.
\medskip

\textit{Claim 2}  $\omega_{f}(y)$ is infinite.

Suppose that $\omega_{f}(y)$ is finite. Again by Lemma \ref{x not pf}, $y\in P(f)$.  Since $\omega_{f}(x)$ is infinite and by Lemma \ref{Deomp}, we can find a periodic arc $I$ such that $\omega_{f}(x) \subset I\cup f(I) \cup \dots f^{s-1}(I)$ and $y\notin I\cup f(I) \cup \dots f^{s-1}(I)$,  where $s$ is the period of $I$. By Lemma \ref{xw}, since $x\in CR(f)$ and $\omega_{f}(x)$ is infinite, $x\in f^{i}(I)$ for some $0\leq i\leq s-1$. Let $\pi$ be a quotient map collapsing each $f^{i}(I), 0\leq i\leq s-1$ to a point. We obtain a factor system $(\pi(X),\pi (f))$ with zero topological entropy on the tree $\pi(X)$ ($\pi$ is monotone) and such that $\pi(x)\neq \pi(y)$ and $\pi(x),\pi(y) \in P(\pi(f))$. Moreover,  $(\pi(x),\pi(y))\in IN(CR(\pi(f)),\pi(f))$. This fact contradicts Claim 1 and thus Claim 2 is proven and both $\omega_{f}(x)$ and $\omega_{f}(y)$ are infinite.
\medskip

Next, as $\omega_{f}(x)$ and $\omega_{f}(y)$ are infinite, by Lemmas \ref{Deomp} and \ref{xw}, the points $x,y$ are contained in some periodic arcs $I_{x}$ and $I_{y}$. Observe that if $I_{x}\cap I_{y}=\emptyset$. Then, by the same arguments of Claim 2, we get a contradiction by collapsing each element of their orbit to a point. Therefore any periodic arc containing $x$ must contain $y$. By Lemma \ref{Deomp}, we may find $a\in \omega_{f}(x)\setminus B(X)$ such that the intersection of all the periodic arcs containing $a$ is reduced to $\{a\}$. So, let $I$ be a periodic arc containing $a$ such that $I\cap B(X)=\emptyset$ and let $m$ be the period of $I$. We have that $\omega_{f}(x)\subset \mathrm{Orb}_{f}(I)$. Thus, by Lemma \ref{xw}, for some $0\leq i\leq m-1,\; x\in f^{i}(I)$ and so $y\in f^{i}(I)$. Recall that $(x,y)\in IN(CR(f),f)$. Thus, by \cite[Theorem 2.9]{tan}, for any $n\geq 0$ we can find $x_{n}\in f^{-n}(x)\cap CR(f)$ and $y_{n}\in f^{-n}(y)\cap CR(f)$, so that $(x_{n},y_{n})\in IN(CR(f),f)$. Clearly for any $n\geq 0,\; x_{n} \neq y_{n}$. Moreover, by Lemma \ref{xw}, each member of the sequences $(x_{n})_{n\geq 0}$ and $(y_{n})_{n\geq 0}$ are contained in $\mathrm{Orb}_{f}(I)$ and they are infinite. Me may pick $n$ large enough so that $x_{n}\in I$ and then $y_{n}\in I$ because otherwise $x$ and $y$ would be in disjoints periodic arcs, which is not possible. Moreover, since $End(I)$ is finite, we may find some $n$ large enough so that $\{x_{n},y_{n}\}\subset I\setminus End(I)$. Without loss of generality, we can assume that $x,y\in I\setminus End(I)$.

Now, let $G=I\cup CR(f)$ and $g=f^{m}|_{G}$. By Lemma \ref{crf}, we have that $CR(f)=CR(g)$ and so $(x,y)\in IN(CR(g),g)$. Without loss of generality we can assume that if $End(I)\cap CR(g)\neq \emptyset$, then $End(I)\cap CR(g) \subset P(g)$. Notice that if we assume that $z\in End(I)\cap CR(g)$ and $z\notin P(g)$, by Lemmas \ref{Deomp} and \ref{xw}, we may find a periodic arc $J\subset I$ with period $q$ containing $z$ and such that $J\cap \{x,y\}=\emptyset$. Then, it is enough to collapse $J$ to a point. Note that it is not restrictive to assume that $End(I)\cap CR(f) \subset Fix(g)$, because we can work with $g^q$ instead of $g$, where $q$ is the period of $z$.

As $End(I)\cap CR(f) \subset Fix(g)$ contains at most two points, we will consider two different cases according to the cardinality of $End(I)\cap CR(g)$. In both cases we will get a contradiction.
\medskip

\textbf{Case.1} $End(I)\cap CR(g)$ contains at most one point.

Assume that either $\{c\}=End(I)\cap CR(f)\neq \emptyset$ and $g(c)=c$ or $End(I)\cap CR(f)= \emptyset$ and $c\in I$ with $g(c)=c$. In both cases,
observe that $Y:=[CR(g)\setminus (I\setminus End(I))]\cup \{c\}\subset G$ is closed. For $z\in Y\setminus I$, by Lemma \ref{xw}, $g(z)\notin I$ and for $z\in [End(I) \cap CR(g)]\cup \{c\}$ we have that $g(z)=z$. Thus, $Y$ is invariant by $g$.
Now, consider the quotient map $\pi$ on $G$ collapsing $Y$ to a point. Clearly $\pi$ induces a factor system $(Z,\phi)$. By Lemma \ref{Factor CRF}, we have that $\pi(CR(g)) \subset CR(\phi)$. Moreover, $\pi(x)\neq \pi(y)$ and, by Lemma \ref{Factor SE(X,f)},  $(\pi(x),\pi(y))\in IN(CR(\phi),\phi)$. Thus, we conclude that $h_{\infty}(\phi|_{CR(\phi)})>0$. Observe that $Z$ is an arc, in fact $\pi(G)=\pi(I)$, and $\pi$ is an homeomorphism on $I$. We conclude that the factor system $(Z,\phi)$ is a continuous map with zero topological entropy on the closed interval $\pi(I)$ and $h_{\infty}(\phi |_{CR(\phi)})>0$, this fact contradicts Lemma \ref{ARC}.
\medskip

\textbf{Case.2} $End(I)\cap CR(g)=\{a,b\}$. Recall that $g(I)\subset I$ and $\{a,b\}\subset Fix(g)$, so let $\psi =g|_{I}$.

First, we claim that  $CR(g)\cap I=CR(\psi )$.
Clearly,  $CR(\psi )\subset CR(g)\cap I$. To prove the converse inclusion, let $z\in CR(g)\cap I$. If $z\in P(g)$ we have finished, then we can assume that $z\notin P(g)$. Thus, $z\in I\setminus End(I)$. If we assume that $z\notin CR(\psi )$, then we would have that $z\notin Q(\psi ,z)$. Since $z\in CR(g)$, then for $n$ large enough any $\frac{1}{n}$-pseudo orbit $\{z,z_{1,n},z_{2,n}, \dots,z_{m_{n},n},z\}$ of $CR(g)$ is not included in $I$.
For any $n\geq 0$, let $1\leq i_{n}\leq m_{n}$ such that $\{z,z_{1,n},z_{2,n}, \dots,z_{i_{n-1},n}\}\subset I$ and $z_{i_{n},n}\notin I$. Without loss of generality, we may assume that the sequences $(z_{i_{n-1},n})_{n\geq 0}$ converges to $z_{1} \in I\cap CR(g)$ and $(z_{i_{n},n})_{n\geq 0}$ converges $z_{2}\in CR(g)\setminus [I\setminus End(I)]$. For any $n\geq 0$, we have that $d(g(z_{i_{n-1},n}),z_{i_{n},n})\leq \frac{1}{n}$. We conclude that $g(z_{1})=z_{2}$ and $g(z_{1})\in I$, thus $g(z_{1})\in I \cap [CR(g)\setminus (I\setminus End(I))]\subset \{a,b\}$. Therefore, $z_{2}\in \{a,b\}$ and since $z_{1}\in CR(g)=CR(f^{m})$, by Lemma \ref{x not pf}, we conclude that $z_{1}=z_{2}\in \{a,b\}$. Without loss of generality, we can assume that $z_{1}=z_{2}=a$. Observe that $(d(g(z_{i_{n-1},n}),a))_{n\geq 0}$ converges to $0$ when $n$ goes to infinity.
For each $n\geq 0$ we may replace $z_{i_{n},n}$ by $a$ so that $\{z,z_{1,n},z_{2,n}, \dots,z_{i_{n-1},n},a\}$ is $\epsilon_{n}$-chain in $I$ from $z$ to $a$, where for each $n\geq 0,\; \epsilon_{n}=d(g(z_{i_{n-1},n}),a)+\frac{1}{n}$. We conclude that $a\in CR(z,\psi )$.
Now, by \cite[Proposition 45]{blockcoppel}, since $a\in CR(z,\psi )$, then either $\psi ^{k}(z)=a$ for some $k\geq 0$ or $a\in Q(\psi,z)$. Recall that $z\in CR(g)\setminus P(g)$. Thus, by Lemma \ref{x not pf}, $\omega_{f}(z)$ is infinite. In particular, $\omega_{\psi}(z)$ is infinite. So, we conclude that $a\in Q(\psi ,z)$. By \cite[Corollary 36]{blockcoppel}, as $a\in Q(\psi ,z)\cap Fix(\psi )$, then we have that $Q(\psi ,z)$ is either a periodic orbit or a finite union of disjoints periodic arcs mapped cyclically by $\psi$. Recall that $Q(\psi ,z)$ contains $\omega_{\psi }(z)$, which is infinite. Moreover, $Q(\psi,z)$ contains a fixed point. Thus, $Q(\psi ,z)$ is an arc $J$ satisfying the following properties: $\psi(J)\subset J$, $\omega_{\psi}(z)\subset J$ and $z\notin J$. Recall that $z\in CR(g)=CR(f)$ and, by Lemma \ref{crf}, $CR(f)=CR(f^{m})$. Moreover, $g=(f^{m})|_{G}$. We conclude that $\omega_{f^{m}}(z) \subset J$ and $z\notin J$. This fact contradicts Lemma \ref{xw}, since $z\in CR(f^{m})$, and the claim is proved.

Let $Y=[CR(g)\setminus (I\setminus End(I))]$. Using the same arguments of Case 1, $Y$ is an invariant-closed subset of $G$. Consider now the quotient map $\pi$ on $G$ collapsing $Y$ to a point. Clearly $\pi$ induce a factor system $(S,\phi)$. Clearly $\pi(x)\neq \pi(y)$ and, by Lemma \ref{Factor SE(X,f)},  $(\pi(x),\pi(y))\in IN(\pi(CR(g)),\phi)$. We conclude that $h_{\infty}(\phi |_{ \pi(CR(g))})>0$. Observe that $\pi(CR(g))=\pi(CR(\psi ))$. In fact, if $z\in CR(g)\setminus (I\setminus End(I))$, we have that $\pi(z)=\pi(a)$ and $a\in CR(\psi )$ and, if not and $z\in I\cap CR(g)$, then by claim 3 we have that $z\in CR(\psi )$. On one hand, $h_{\infty}(\phi |_{ \pi(CR(g))})>0$, and on other hand, by Lemma \ref{ARC}, we have that $h_{\infty}(\psi |_{CR(\psi )})=0$. Then $h_{\infty}(\phi |_{\pi(CR(\psi ))})=0$, which is a contradiction because $\pi(CR(g))=\pi(CR(\psi ))$.
\qed

\begin{rem}\label{rm1}
\rm{
For one-dimensional dynamics, the notions of topological sequence entropy and chaos in the sense of Li and Yorke are deeply connected. If $X$ is a compact interval, a circle or a tree and $f:X\rightarrow X$ is continuous with zero topological entropy, then $h_{\infty}(f)>0$ if and only if $f$ is chaotic in the sense of Li and Yorke (see \cite{smifra,hric2,df}). In addition, by \cite{canovas,canovas3,tan}, the supremum of the topological sequence entropies of chaotic maps $f$ with zero topological entropy is $h_{\infty}(f)=\log2$. The result of this paper, jointly with \cite{canovas3}, show that, for continuous self maps defined on trees, the chaos is not located in their chain recurrent set. A similar result was proved for circle maps in \cite{kuang} by replacing the chain recurrent set with the non-wandering set. It is unclear whether this result is possible for circle and, in general, for continuous graph maps when considering the chain recurrent set.
}
\end{rem}

\begin{rem}
\rm{
If $\mu$ is an invariant measure of a continuous tree map $f$, (recall that $\mu$ is invariant if $\mu (f^{-1}(Y))=\mu (Y)$ for any Borel set $Y\subset X$), then the supremum of metric sequence entropies of $f$ (see \cite{kusni} for the definition of metric sequence entropy) is zero. This fact is because the support of any invariant measure is contained in the non-wandering set of $f$ (cf. \cite[Chapter 6]{walters}), and topological sequence entropy is an upper bound for metric sequence entropy \cite{goodman}.  }
\end{rem}

\section{A counterexample for dendrite maps}
To conclude this paper, we give an example of a dendrite map with zero topological entropy such that $\Omega(f)=P(f)$ and $h_{\infty}(f|_{\Omega(f)})>0$. The idea as in \cite{IA} is to consider a pointwise periodic homeomorphism $g$ on a given compact set $Y$ such that $h_{\infty}(g)>0$ and then extend the homeomorphism to a continuous self-map on the whole dendrite. Note that in \cite{IA} the homeomorphism $g$ was extended to a pointwise periodic homeomorphism, but the resulting space was a dendroid. This fact cannot happen in the case of dendrite homeomorphism since any dendrite homeomorphism, and in particular, a tree homeomorphism, has zero sequence topological entropy (see \cite{es0}).

Let $I=[0,1]\times \{0\}$. For $n\geq 1$, let
$$
Y_{n}:= \left\{ \left( x_{i,n},\frac{1}{n} \right), \left( x_{i,n},\frac{-1}{n} \right): i\in \{ 1,2,...,2^{n-1} \} \right\} ,
$$
where $x_{i,n}:= \frac{2i-1}{2^n}$ for $n\geq 1$ and $ 1\leq i\leq 2^{n-1}$.
Observe that
$$
\lim_{n\to +\infty} Y_{n}=I
$$
and $d(x_{i,n},x_{i+1,n})\leq \frac{1}{2^{n-1}}$ for $1\leq i\leq 2^{n-1}-1$.

Let $Y:=I\cup \left(\bigcup_{n\geq 1}Y_{n} \right)$, which clearly is a compact subset of $\mathbb{R}^2$. We define $g:Y\to Y$ by
\begin{enumerate}
\item[(i)] $g(y)=y$, if $y\in I$.
\item[(ii)] For any $n\geq 1$ the subset $Y_{n}$ is a periodic orbit such that:
\medskip
\begin{itemize}
\item $g((x_{i,n},\frac{1}{n}))=(x_{i+1,n},\frac{1}{n})$ if  $1\leq i\leq 2^{n-1}-1$.
\item $g((x_{2^{n-1},n},\frac{1}{n}))=(x_{2^{n-1},n},\frac{-1}{n})$.
\item $g(x_{i,n},\frac{-1}{n})=(x_{i-1,n},\frac{-1}{n})$  if  $2\leq i\leq 2^{n-1}$.
\item $g((x_{1,n},\frac{-1}{n}))=(x_{1,n},\frac{1}{n})$.
\end{itemize}
\end{enumerate}
\medskip

Clearly, $g$ is a pointwise periodic homeomorphism and thus $(Y,g)$ is distal. We claim that $h_{\infty}(g)>0$. Otherwise if $h_{\infty}(g)=0$, then $(Y,g)$ will be a distal null system. Therefore by \cite[Theorem 3.2]{null}, the map $g$ would be Lyapunov stable, which is not the case since the powers of $g$ are not equicontinuous on any point of $I$.

\begin{figure}[!h]
 \centering
  \includegraphics[width=12cm, height=7.5cm]{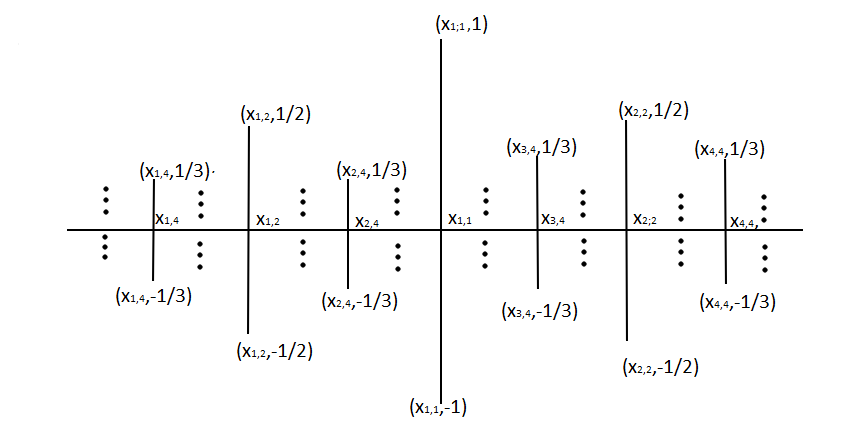}
  \caption{The dendrite $D$}
\label{ent}
\end{figure}

Now, we extend the map $g$ on a dendrite as follows. For any $n\geq 1$ and $i\in \{ 1,2,...,2^{n-1} \}$, let $I_{n,i}$ be the arc joining $(x_{n,i},\frac{1}{n})$ to $(x_{n,i},0)$ and let $J_{n,i}$ be the arc joining $(x_{n,i},\frac{-1}{n})$ to $(x_{n,i},0)$. Then, $D=I \cup \big(\displaystyle\bigcup_{n\geq 1, 1\leq i\leq 2^{n-1}}I_{n,i}\cup J_{n,i}\big)$ is a dendrite (see Figure \ref{ent}).
We extend $g$ into a map of $D$ as follows:\\
\begin{enumerate}
\item[(i)] For $1\leq i\leq 2^{n-1}-1,\; f$ maps $I_{n,i}$ to the arc $[(x_{n,i},0),(x_{n,i+1},0)]\cup I_{n,i+1}$ in an affine way.
\item[(ii)] $f$ maps $I_{n,2^{n-1}}$ into $J_{n,2^{n-1}}$ in an affine way and maps $J_{1,n}$ to $I_{1,n}$ in an affine way.
\item[(iii)] For $2\leq i\leq 2^{n-1},\; f$ maps $J_{n,i}$ to the arc $[(x_{n,i},0),(x_{n,i-1},0)]\cup I_{n,i-1}$ in an affine way.
\item[(iv)] $f|_{A}=g$.
\end{enumerate}

\medskip

The map $f$ is continuous on $D$. Moreover, any point of $D$ is periodic or eventually periodic, and thus $h(f)=0$. In addition, $\Omega(f)=P(f)=Y$ and $h_{\infty}(f)\geq h_{\infty}(f|_{Y})=h_{\infty}(g)>0$. It is also worth mentioning that, although $h_{\infty}(f)>0$, the map $f$ is not chaotic in the sense of Li and Yorke. It does not have any Li-Yorke pairs since any point is periodic or eventually periodic. This fact cannot happen in the case of  tree maps, because positive sequence entropy characterizes the existence of chaos in the sense of Li and Yorke for continuous tree maps (see Remark \ref{rm1}).

\section*{Acknowledgements}

Aymen Daghar was supported by the research unit: ``Dynamical systems and their applications'', (UR17ES21), Ministry of Higher Education and Scientific Research, Faculty of Science of Bizerte, Bizerte, Tunisia.
This work was made during the Aymen Daghar's visit at Universidad Polit\'{e}cnica de Cartagena under the Erasmus program KA107. The support of this university is also gratefully acknowledged.

J.S. C\'{a}novas was supported by the project MTM 2017-84079-P Agencia
Estatal de Investigaci\'{o}n (AEI) y Fondo Europeo de Desarrollo Regional
(FEDER).


\bigskip

\bibliographystyle{amsplain}

\end{document}